\documentclass[12pt,a4paper,oneside]{article}
\newif\ifRUS\RUSfalse\newif\ifENG\ENGtrue
\newif\ifDRAFT\DRAFTfalse
\usepackage{amsmath}    
\usepackage{amsthm}     
\usepackage{amsfonts,amssymb}
\usepackage{cmap}  
\ifENG
\usepackage[T2A]{fontenc}    
\usepackage[russian,english]{babel}
\else
\usepackage[T1,T2A]{fontenc}    
\usepackage[english,russian]{babel}      
\fi
\usepackage{graphicx}
\newcommand{\Skip}[1]{}
\usepackage{fullpage}
\usepackage{ifthen}
\usepackage{textcomp} 
\usepackage{xcolor}
\usepackage{caption}
\ifRUS
\usepackage{indentfirst}
\fi
\usepackage[medium,center,small]{titlesec}
\titleformat
{\chapter}                
[display]                 
{\bfseries\Large\itshape} 
{Story No. \ \thechapter} 
{0.5ex}                   
{
    \rule{\textwidth}{1pt}
    \vspace{1ex}
    \centering
} 
[
\vspace{-0.5ex}%
\rule{\textwidth}{0.3pt}
] 
\titlelabel{\thetitle.~}
\newenvironment{myitemize}
{\begin{itemize}
    \setlength{\itemsep}{0pt}
    \setlength{\parskip}{0pt}
    \setlength{\parsep}{0pt}
}
{\end{itemize}
}
\newsavebox{\tmpbox}

\newlength{\tmplength}

\newcommand{\textwh}{\text{\ifENG{where}\else{где}\fi}}  
\newcommand{\textif}[1][]{\text{#1\ifENG{if}\else{если}\fi#1}}

\newcommand{\ie}{\ifENG{i.\,e.}\else{т.\,е.}\fi}
\newcommand{\quo}[1]{\glqq#1\grqq}
\newcommand{\tire}{~\textbf{---}\ }

\newcommand{\MP}[1]{\marginpar{\color{magenta}\small#1}}
\ifDRAFT
\newcommand{\ShowL}{\reversemarginpar\MP{~\the\inputlineno}\normalmarginpar}
\else
\newcommand{\ShowL}{}
\fi

\DeclareMathOperator{\sinc}{sinc}
\ifRUS

\let\arctan\relax
\DeclareMathOperator{\arctan}{arctg}

\let\th\relax%
\DeclareMathOperator{\th}{th}

\else  
\fi
\newcommand{\abs}[1]{\left\lvert#1\right\rvert}
\newcommand{\Deg}[1]{{\ifmmode{#1}^\circ\else{#1}\textdegree\fi}}
\DeclareMathOperator{\sgn}{sgn}
\newcommand{\Exp}[1]{\mathrm{e}^{#1}}

\newcommand{\iu}{{\mathrm{i}}}  
\newcommand{\Brack}[1]{\left[\,#1\,\right]}
\newcommand{\Brace}[1]{\left\{#1\right\}}
\newcommand{\Skobki}[1]{\left(#1\right)}
\newcommand{\where}[1][\textwh]{\quad\text{#1}\quad}
\newcommand{\ALI}[1]{\aligned #1\endaligned} 
\newcommand{\Equa}[2]{%
\ifDRAFT\marginpar{\vspace*{1\baselineskip}{\color{brown}\small~#1}}\fi%
\begin{equation}#2\label{#1}\end{equation}%
}
\newcommand{\equa}[1]{\[ #1 \]}
\ifRUS
\renewcommand{\le}{\leqslant}
\renewcommand{\ge}{\geqslant}
\fi
\newcommand{\Eqref}[1]{Eq.\,\ref{#1}}

\newcommand{\eqpipi}[2]{#1\equiv#2\!\!\pmod{2\pi}} 
\newcommand{\vp}{\varphi}

\newcommand*\Bbar[1]{
  \hbox{%
    \vbox{%
      \hrule height 0.5pt 
      \kern0.5ex
      \hbox{%
        \kern-0.1em
        \ensuremath{#1}%
        \kern-0.1em
      }%
    }%
  }%
}
\newcommand{\St}[2][]{{#2}_{#1}^{\star}}
\newcommand{\So}{\quad\Longrightarrow\quad} 
\newcommand{\Vect}[1]{\protect\overrightarrow{#1}} 

\newcommand{\half}[1]{\dfrac{#1}2}
\newcommand{\acum}{\widetilde\alpha}
\newcommand{\bcum}{\widetilde\beta}
\newcommand{\Kl}[2][]{{\ifmmode{\cal K}^{#1}_{#2}\else${\cal K}^{#1}_{#2}$\fi}}
\newcommand{\Kls}[1]{{\ifmmode{\cal K}^{\star}_{#1}\else${\cal K}^{\star}_{#1}$\fi}}

\newcommand{\Sth}[1][]{\sin^{#1}\theta}
\newcommand{\Cth}[1][]{\cos^{#1}\theta}
\ifDRAFT
\usepackage{fancyhdr}
\def\mydate{\leavevmode\hbox{\the\year-\twodigits\month-\twodigits\day}}
\def\twodigits#1{\ifnum#1<10 0\fi\the#1}
\fi
\newcommand{\Figref}[1]{\ref{fig:#1}}
\newcommand{\Reffig}[2][]{\ifRUS рис.\else fig.\fi\,\Figref{#2}{#1}}
\newcommand{\RefFig}[2][]{\ifRUS Рис.\else Fig.\fi\,\Figref{#2}{#1}}
\newcommand{\Infigw}[2]{
\includegraphics[width=#1]{#2.eps}}
\newcommand{\Infig}[3]{\Infigw{#1}{#2}\caption{{\small #3}}\label{fig:#2}}
\newcommand{\Eqfig}[2]{
\ifmmode%
\parbox[c]{#1}{\Infigw{#1}{#2}}%
\else
\ifDRAFT\marginpar{\vspace*{1\baselineskip}{\color{brown}\small~#2}}\fi%
\equa{\parbox[c]{#1}{\Infigw{#1}{#2}}}%
\fi%
}
\newcommand{\Pfig}[3]{
\centering%
\ifDRAFT%
\settowidth{\tmplength}{#2}\addtolength{\tmplength}{\textwidth}%
\makebox[\tmplength][r]{{\color{magenta}\small#2}}\par\vspace{-\baselineskip}%
\fi%
\Infig{#1}{#2}{#3}}
\begin{document}
%
%
\setlength{\abovedisplayskip}{6pt plus 1pt minus 1pt}
\setlength{\abovedisplayshortskip}{6pt plus 2pt}
\setlength{\belowdisplayskip}{4pt plus 1pt minus 1pt}
\setlength{\belowdisplayshortskip}{4pt plus 2.5pt minus 2.0pt}
%
\ifRUS
\else
\selectlanguage{english}%
\author{Alexey Kurnosenko (\,\texttt{alexeykurnosenko@gmail.com}\,)} 
\title{Construction of a spiral with given boundary conditions by inversion of the involute of a circle}
\fi
\maketitle{}
%
\vspace{-2\baselineskip}%
\begin{abstract}
\ifRUS
\else
To construct a curve with a monotonic curvature (spiral),
and given tangents and curvatures at the ends,
the author proposed the following method.
From given boundary conditions, the values of two inverse invariants are determined.
Then, on some base spiral (initially, a logarithmic spiral was chosen),
an arc with the same invariant values is sought for.
A linear-fractional map of the found arc solves the problem.
It seems that choosing the involute of a circle as the base spiral
yields the simplest solution, which we present here.
\fi
\end{abstract}

\ShowL%
\ifRUS
\else
The problem of constructing a curve with a monotonic change in curvature (a spiral)
in Computer Aided Design (CAD) applications is commonly referred to as
"two-point G2 Hermite interpolation with spirals."
The $G^1$ version of the problem requires construction with specified tangents at the curve's ends,
while the $G^2$ version specifies also boundary curvatures.
Additional conditions may be: the requirement of curvature continuity,
a solution in terms of curves, accepted in CAD applications (e.g., Bézier curves),
constraints on the curve's rotation.
\fi

\ifRUS
\else
The simplest solution of the $G^2$ problem is to construct a triarc,
\ie{} a smooth curve,
consisting of three circular arcs
\cite{TriarcsMW,InvInv}.
Other solutions are proposed in
\fi
\cite{MeekPlanar,FreyConic,DietzRational,GoodmanInvolute}.

\ShowL%
\ifRUS
\else
In \cite{PomiG2} the author proposed a method for constructing an analytical spiral
with any boun\-dary conditions, under which monotonicity of curvature is possible,
including spirals with an inflection point, and spirals, curling around end points.
Construction of the spiral is performed as follows.
\fi
\begin{myitemize}
\item
\ifRUS
\else
From boundary conditions,
the values of inverse invariants $Q$ and $\omega$ \eqref{DefQomg} are determined.
\fi
\item
\ifRUS
\else
On a certain base spiral an arc is sought for,
that realizes the required invariant values.
\fi
\item
\ifRUS
\else
A linear-fractional map of the found arc solves the problem,
preserving the monoto\-ni\-city of curvature and analyticity of the curve.
\fi
\end{myitemize}

\ifRUS
\else
Different choices of the base spiral yield a variety of solutions,
allowing the designer to select a particular arc length,
or a preferable function profile $k(s)$ in the natural equation of the curve.
A transcendental curve, typically obtained by this method,
can easily be transformed into a curve with piecewise constant curvature.
\fi

\ShowL%
\ifRUS
\else
The logarithmic spiral was chosen as the base spiral in \cite{PomiG2}.
The solution is a double logarithmic spiral, that matching any boundary conditions.
\fi

\ifRUS
\else
In \cite{AKhyperb} the base spiral was an arc of the hyperbola between two vertices.
This method is applicable for relatively short spirals.
The result is a rational 4th order curve.
\fi

\ifRUS
\else
The choice of the involute of the circle also allows to satisfy any boundary conditions,
and seems to provide the simplest solution, which we present here.
\fi
\smallskip

\ShowL%
\ifRUS
\else
To construct the involute of a circle of radius~$R$,
we use the model of a thread, unwinding from a circle.
The coordinate origin is placed at the center of the circle.
The parametric equation of the involute is $[x(t),\,y(t)]$,
where~$t$ is the
polar angle of the point of contact of the thread and the circle.
$L=Rt$ is the length of unwound thread,
$s(t)$ is the length of the arc of involute,
$\tau(t)$ and $k(t)$ are the slope of the tangent and the curvature;
$\rho$~and~$\vp$ are polar coordinates of the curve:
\fi
\equa{%
  \ALI{
     &x(t)=R\cos t+L\cos(t{-}\pi/2)=R(\cos t+ t\sin t),\\
     &y(t)=R\sin t+L\sin(t{-}\pi/2)=R(\sin t- t\cos t);\\
    &s(t)=\frac12 R t^2,\quad \tau(t) = t,\quad k( t)=\dfrac{1}{R t},\;
     k(s) = \dfrac{1}{\sqrt{2Rs}}.\\
    &\rho( t)=R\sqrt{ t^2+1},\quad\varphi( t)= t-\arctan t.
  }
  \;
  \Eqfig{200bp}{EvolvDef}
}
\ifRUS
\else
The curvature of the involute is positive and decreases from $+\infty$ to $0$.
Note that the circles of curvature at $t=t_1$ and $t_2=t_1+2n\pi$ are concentric.
\fi

\ShowL%
\ifRUS
\else
Let us clarify what is meant by boundary conditions.
We will represent the sought for curve,
subtented by the chord $AB$ of length $\abs{AB}=2c$,
in the coordinate system of the chord,
in which the origin is in the middle of the chord,
and the direction of the vector~$\Vect{AB}$
coincides with the positive direction of the $X$-axis.
The coordinates of the starting and ending points are:
${A=(-c,0)}$, ${B=(c,0)}$.
The slopes of the tangents at the points~$A$ and~$B$ and are denoted by  $\alpha$~and~$\beta$,
boundary curvatures are $k_1$ and~$k_2$.
The type of monotonicity of curvature (increasing, decreasing) is denoted by
${M=\sgn(k_2-k_1)=\pm1}$.
Angles are determined within the following limits:
\fi
\Equa{A0B0}{
    \aligned
    &\textif\quad M=+1{:}\quad \alpha,\,\beta\in(-\pi;\,\pi];\\
    &\textif\quad M=-1{:}\quad \alpha,\,\beta\in[-\pi;\,\pi);
    \endaligned\qquad
    \acum=\alpha+2Mn_1\pi,\quad \bcum=\beta+2Mn_2\pi,\quad n_{1,2}\in\mathbb{N}_0.
}
\ShowL%
\ifRUS
\else
The tilde marks the cumulative versions of the angles, $\acum$ и $\bcum$,
reflecting the degree of curling
of spirals around the ends: $n_1$ and $n_2$
are the number of intersections of the internal points of the spirals
with rays, complementing the chord on the left and right to an infinite straight line.
\RefFig{ThreeExa} 
illustrates the definition of cumulative angles.
The difference $\bcum-\acum$ is equal to a full rotation of the spiral
(the integral of the curvature $k(s)$ over the arc length~$s$).
The half-intervals, specified in \eqref{A0B0} for $\alpha$ and $\beta$,
are such that for a short spiral,
i.e. an arc that does not intersect the complement of the chord,
equalities $\acum=\alpha$ and $\bcum=\beta$ are fulfilled.
\fi
\ifRUS
\else
Note that these quantities are quite natural:
if a small arc of a spiral, for which~$\alpha$ and~$\beta$ are close to zero,
is extended to a full arc,
maintaining continuity when passing through values that are multiples of~$\pm\pi$, 
ultimately yields values  $\acum$~and~$\bcum$  \cite[Corollary\,5.1]{PomiMain}.
An arc of involute
curls only around the start point: $n_2=0$
(arc $A_3B_3$ in \RefFig{Chords}).
\fi

\begin{figure}[t]
\Pfig{\textwidth}{ThreeExa}
{
\ifRUS
\else
Three spirals with the same values of $\Brace{c,\alpha,\beta,k_1,k_2}$,
but different $\acum,\bcum$
\fi
}
\end{figure}

\Skip{ Equa{K1K2cc}%
      \Kl{1} = K(-c,0,\alpha,k_1)  \qquad\mbox{и}\qquad
      \Kl{2} = K(c,0,\beta,k_2).
}

\ShowL%
\ifRUS
\else
Define now two invariants with respect to M\"obius maps:
\fi
\Equa{DefQomg}{%
    \omega=\frac{\acum+\bcum}2;\quad 
    Q = 
     (k_1c+\sin\alpha)(k_2c-\sin\beta)+\sin^2\frac{\alpha{+}\beta}{2}
}
\ifRUS
\else
(under inversion $M$ and $\omega$ change sign).
$Q$ and $\omega$  depend on the angles between oriented circles.
The quantity $2\omega$ is the angular width of the lens, bounding the short spiral arc
\hbox{\cite{PomiShort,PomiLong}}.
The relationship of~$Q$ with the inverse distance between circles is commented in
\cite{InvInv}.
The necessary conditions for the existence of a spiral with given boundary conditions are:
$\sgn\omega=M$ \cite[th.\,1]{PomiLong} and $Q\le0$ \cite[th.\,2]{PomiMain};
the equality $Q=0$ arises only if the boundary circles of curvature are tangent
(the curve is a biarc).
\fi
\smallskip

\ShowL%
\ifRUS
\else
Next we set $R=1$, and work with the involute with increasing curvature (from $-\infty$ to~0),
for which we symmetrically reflect the curve relative to the abscissa axis:
\fi
\Equa{IncrCrv}{
  y(t)\to{-y(t)},\quad \tau(t)\to{-\tau(t)},\quad k(t)\to{-k(t)}.
}
%
%
\ifRUS
\else
We are looking for an arc $t\in[t_1,\,t_2]$,
on which the required values $\omega$ and $Q$ are realized.
Let
\fi
\equa{
  t_1=t_0-\theta>0,\quad  t_2=t_0+\theta>t_1.   
}
\ifRUS
\else
Let's determine the length of the chord $2c=\abs{AB}$
and the direction~$\mu$ of vector~$\Vect{AB}$:
\fi
\equa{
  c = \sqrt{(\theta\cos\theta-\sin\theta)^2+t_0^2 \sin^2\theta},\quad
  \aligned
    &c\,\cos\mu=
     \theta\sin t_0\cos\theta +\Skobki{t_0\cos t_0-\sin t_0}\sin\theta,\\
    &c\,\sin\mu=
     \theta\cos t_0\cos\theta -\Skobki{t_0\sin t_0+\cos t_0}\sin\theta.
  \endaligned
}
\ifRUS
\else
Substituting into \eqref{DefQomg}
\fi
\equa{
  \eqpipi{\alpha}{\Brack{\tau(t_1)-\mu}}=-t_1-\mu,\qquad
  \eqpipi{\beta }{\Brack{\tau(t_2)-\mu}}=-t_2-\mu,
}
\ifRUS
\else
we find:
\fi
\Equa{EQt0}{
    Q=\frac{\theta^2-\sin^2\theta}{\theta^2-t_0^2}
    \So
    t_0(\theta)=\sqrt{\frac{\theta^2(1-Q)-\sin^2\theta}{-Q}}.
}
\ShowL%
\ifRUS
\else
The required $Q$ is realized on a family of chords with the found dependence $t_0(\theta)$.
In \RefFig{Chords} such family is shown.
Three chords are highlighted, and shown, together with the arc of the involute,
on the right fragment of the figure in the chord coordinate system.
Let's find a chord in the family that gives the required value of
\fi
\ $\eqpipi{\omega}{\Brack{-\mu(\theta)-t_0(\theta)}}$:
\equa{
 \ALI{
  &
   \cos\omega=\frac{t_0(\theta)\sin\theta}c,\quad
    \sin\omega=\frac{\sin\theta-\theta\cos\theta}c,\quad
    \tg\omega=\frac{\sin\theta-\theta\cos\theta}{t_0(\theta)\sin\theta};\\
  &\tg\Skobki{\omega-\theta+\frac{\pi}2}=\sin\theta\cdot
     \frac{(t_0(\theta)-\theta)\cos\theta+\sin\theta}
          {(t_0(\theta)-\theta)\sin^2\theta+\theta-\cos\theta\sin\theta}.
 }
}
\begin{figure}[t]
\Pfig{1\textwidth}{Chords}%
{
\ifRUS
\else%
The family of chords of the involute \eqref{IncrCrv}, on which  $Q=-0.04$;
on the right three selected chords are shown in their own coordinate system
\fi%
}
\end{figure}
\ShowL%
\ifRUS
\else
Since the denominator on the right-hand side is positive $\Brack{t_0(\theta)>\theta}$,
the function is con\-ti\-nuous.
The final equation for~$\theta$:
 \fi
\Equa{eqF}{
  \ALI{%
   \omega(\theta)&{}=
   \theta-\half\pi+
   \arctg\frac{\Sth\,\Brack{t_0(\theta)\Cth-\theta\Cth+\Sth}}%
              {t_0(\theta)\sin^2\theta+\Cth\,(\theta\cos\theta-\sin\theta)}.\\
  }
}
\ifRUS
\else
On the right-hand side there should be a term, that is a multiple of $2\pi$,
but it is equalized to zero,
since $\omega(0) = 0$ and the function is continuous:
\fi
\equa{
    \lim_{\theta\to0}\arctg\Brack{\frac{\ldots}{\cdots}}=
     \lim_{\theta\to0}\,\arctg\Brack{\frac32\theta^{-1}+\frac{8Q-5}{40Q}\theta+O(\theta^3)}=\half\pi\,.
}
\ShowL%
\ifRUS
\else
The graphs of the function $\omega(\theta)$ for different~$Q$ are shown in \RefFig{ThetaBaseQ}.
The function is monotone: its derivative, written with an accuracy up to a positive factor,
is equal to
\fi
\equa{
  \sinc^2\theta\Skobki{\Sth[2]+2\sinc\theta\Cth -3-\abs{Q}}+1+\abs{Q}\quad
  \Brack{\sinc\theta\equiv\frac{\sin\theta}{\theta}}.
}

\begin{figure}[t]
\Pfig{1.\textwidth}{ThetaBaseQ}%
{
\ifRUS
\else
Graphs $\omega(\theta)$ \eqref{eqF}; the dotted line is the graph for
\fi
$Q={-\infty}$~ $\Brack{t_0(\th) = \th}$
}
\end{figure}

\ifRUS
\else
Numerical solution of \Eqref{eqF}, for example, by the bisection method, is easy.
The roots of the equation~$\tg\theta=\theta$ are also the roots of equation~\eqref{eqF}:
\fi
\equa{
    \theta_0=0,\quad
    \theta_1\approx 1.430\pi,\quad 
    \theta_2\approx 2.459\pi,\quad
    \theta_3\approx 3.471\pi,\quad\ldots,\quad
    \theta_n \approx\frac{(2n+1)\pi}2-\frac2{(2n+1)\pi}. 
}
\ifRUS
\else
So, the search interval can be reduced to
\fi
$\Brack{\theta_i;\,\theta_{i+1}}$.

\Skip{
Тогда, если $\omega=n\pi$, то $\Theta_n$\tire решение уравнения~\eqref{eqF};
если $\omega\in[n\pi; (n+1)\pi]$, то решения следует искать на отрезке
 $\theta\in[\Theta_n;\Theta_{n+1}]$  ($I_0\in[0; 3\pi/2]$, $I_n\in[(n+0.4)\pi; (n+1.5)\pi]$).
\equa{
  \aligned
    &\omega\in(0;\,\pi/2]&&                          && \So &&\theta\in(2\omega;\,\pi]; \\
    &\omega\in\Skobki{n\pi-\pi/2;\,n\pi+\pi/2},\;&& n>0&& \So &&\theta\in\Skobki{n\pi;\,(n+1)\pi};\\
    &\omega = \pi/2+n\pi,&&                     n\ge0&& \So &&\theta=(n+1)\pi.
  \endaligned
}
\verb!fsolve(tan(x*Pi)-x*Pi, x, 1..10,maxsols=5);!
}

\begin{figure}[t]
\Pfig{.92\textwidth}{EvCornu}%
{
\ifRUS
\else
The boundary conditions are borrowed from an arc of the Cornu spiral (dotted line).
The curve, labelled as {E}, is derived from the involute of a circle.
The curves, labelled {L} (Log. spiral) and {H} (Hyperbola),
are derived from other base spirals, shown at the top.
\fi
}
\end{figure}

\ifRUS
\else
For values $Q,\omega$, calculated from the given boundary conditions,
we define~$\theta$~\eqref{eqF} and~$t_0$~\eqref{EQt0},
find the arc $\Brack{t_1;\,t_2}$ of the involute.
By rotation, translation and scaling
(dividing coordinates by~$c$ and multiplying curvatures by~$c$)
we bring it into a coordinate system
with the unit chord  ($c=1$).
The boundary conditions of the resulting arc are denoted as
$\St{\alpha}$, $\St{\beta}$, $\St[1]{k}$, $\St[2]{k}$.
The M\"obius map of the obtained arc $z(t) = x(t) + \iu\,y(t)$ of the involute
to the desired one, also reduced to the unit chord, has the form:
\fi
\Equa{Moebius}{
\hspace*{-5mm}%
        \frac{z(t)+z_0}{1+z_0\, z(t)},\where
    z_0=\frac{\rho\Exp{\iu\lambda}-1}{\rho\Exp{\iu\lambda}+1},\quad
      \lambda=\alpha^\star-\alpha = \beta-\beta^\star\,,\quad
      \rho=\frac{k_1 + \sin\alpha}{\St[1]k + \sin\St{\alpha}}=
        \frac{\St[2]k - \sin\St\beta}{k_2 - \sin\beta}.
}

\ifRUS
\else
Three examples are show in Figs \Figref{EvCornu}, \Figref{EvConc}, \Figref{EvTrx}.
\fi

\ShowL%
\ifRUS
\else
In \RefFig{EvCornu},
the boundary conditions are borrowed from an arc of the Cornu spiral,
shown in the lower fragment of the figure by the dotted line.
The dashed circles are the boundary circles of curvature.
The curve, marked with the symbol E, is obtained by the present method.
Curves, marked as L (Log. spiral) and H (Hyperbola),
were obtained in a similar way with other base spirals \cite{PomiG2,AKhyperb}.
The middle fragment shows the curvature graphs $k(s)$ of the obtained interpolants.
The upper fragment of the figure illustrates three arcs of base spirals,
from which, by transformation~\eqref{Moebius}, three solutions were obtained.
\fi

\ShowL%
\ifRUS
\else
\clearpage
In \RefFig{EvConc} the boundary conditions are a pair of concentric circles.
In the first configuration the solutions E,L,H are practically indistinguishable
(only E is shown).
In the second one there is no solution H (the curve is too long, $\omega>\pi/2$). 
\fi

\ShowL%
\ifRUS
\else
In \RefFig{EvTrx}, the boundary conditions are taken from an arc of the tractrix of a circle,
namely, the polar (or spiral) tractrix:
the length of the leash of this tractrix is equal to the radius of the circle.
This is the curve with a natural equation and a natural parameterization
\fi
\equa{
  \ALI{%
    &k(s)= \frac{1-2\Exp{-s/T}}{T\sqrt{1-\Brack{1-2\Exp{-s/T}}^2}};\\
    &\tau(s)=\sqrt{\Exp{s/T}-1}-\arccos\Skobki{2\Exp{-s/T}-1}+\pi,\\
    &x(s)=2T\Exp{-s/T}\Brack{\psi(s)\sin\psi(s)+\cos\psi(s)},\\
    &y(s)=2T\Exp{-s/T}\Brack{\sin\psi(s)-\psi(s)\cos\psi(s)},\\
    &\textwh\quad \psi(s)=\sqrt{\Exp{s/T}-1}.
  }\quad
  \Eqfig{180bp}{TrxDef}
}
\ifRUS
\else
$P$ is the inflection point of the tractrix.
Since this curve is the inverse of the involute,
the solution~E exactly coincides with the original tractrix,
as can be seen in the figure.
\fi

\begin{figure}[t]
\Pfig{.8\textwidth}{EvConc}%
{
\ifRUS
\else
The boundary circles of curvature are concentric.
\fi
}
\end{figure}

\begin{figure}[h]
\Pfig{.8\textwidth}{EvTrx}%
{
\ifRUS
\else
Boundary conditions are taken from an arc of polar tractrix (dotted curve)
\fi
}
\end{figure}

\end{document}